\documentclass{amsart}

\newtheorem{theorem}{Theorem}[section]

\theoremstyle{definition}

\theoremstyle{remark}
\newtheorem{remark}{Remark}[section]

\numberwithin{equation}{section}

\newcommand{\abs}[1]{\left\lvert#1\right\rvert}

\usepackage[shortlabels]{enumitem}

\newcommand{\eulerianNum}[2]{\left\langle \genfrac{}{}{0pt}{}{#1}{#2} \right\rangle}

\usepackage{graphicx}

\usepackage{orcidlink}

\begin{document}

\title[Functional equation for certain LC-functions]{Functional equation for LC-functions with even or odd modulator}

\author[L. Lamgouni]{Lahcen Lamgouni \orcidlink{0000-0002-3927-0353}}

\curraddr{Ministry of National Education Preschool and Sports, Errachidia, Morocco}
\email{lahcen.lamgouni@taalim.ma}





\begin{abstract}
In a recent work, we introduced \textit{LC-functions} $L(s,f)$, associated to a certain real-analytic function $f$ at $0$, extending the concept of the Hurwitz zeta function and its formula. In this paper, we establish the existence of a functional equation for a specific class of LC-functions. More precisely, we demonstrate that if the function $p_f(t):=f(t)(e^t-1)/t$, called the \textit{modulator} of $L(s,f)$, exhibits even or odd symmetry, the \textit{LC-function formula}---a generalization of the Hurwitz formula---naturally simplifies to a functional equation analogous to that of the Dirichlet L-function $L(s,\chi)$, associated to a primitive character $\chi$ of inherent parity. Furthermore, using this equation, we derive a general formula for the values of these LC-functions at even or odd positive integers, depending on whether the modulator $p_f$ is even or odd, respectively. Two illustrative examples of the functional equation are provided for distinct parity of modulators.
\end{abstract}

\maketitle

\section{Introduction}\label{s:Int}

The Hurwitz zeta function $\zeta(s,a)$ is defined for complex variables $s$ with $\Re(s)>1$ and $a\in\mathbb{C}\setminus\mathbb{Z}_{\leq0}$ by the series $\zeta(s, a)=\sum_{n=0}^{+\infty}(n+a)^{-s}$. This function can be analytically extended to all $s\in\mathbb{C}$ except at $s=1$, and satisfies the well-known Hurwitz formula, for $0<a\leq 1$ and $\Re(s)>1$:
\begin{equation}
\label{e:Hur_For}
\zeta(1-s,a)=\frac{\Gamma(s)}{(2\pi)^{s}}\left(e^{-\frac{i\pi s}{2}}F(a, s)+e^{\frac{i\pi s}{2}}F(-a,s) \right).
\end{equation}
Here $\Gamma(s)$ is the Gamma function and $F(a,s)$ is the periodic zeta function defined for $\Re(s)>1$ by $F(a,s)=\sum_{n=1}^{+\infty}e^{2in\pi a}n^{-s}$; for more details, see, e.g., \S 12 of \cite{MR0434929}. Note that $F(a,s)$ for $0<a<1$ has analytic continuation to the whole complex plane (see, e.g., \cite{MR43843, MR1811017}).

The Dirichlet L-function, $L(s,\chi)$, on the other hand, is defined for a Dirichlet character $\chi$ modulo $q$, where $q\geq1$, by the series 
\[
L(s,\chi)=\sum_{n=1}^{+\infty}\chi(n)n^{-s},
\]
which converges for $\Re(s)>1$. $L(s,\chi)$ extends to an entire function for a non-principal character $\chi$ (i.e. $\sum_{k=1}^{q}\chi(k)=0$), otherwise to a meromorphic function with a simple pole at $s=1$, with residue $\sum_{k=1}^{q}\chi(k)/q=\varphi(q)/q$ (see, e.g., \cite[Theorem 12.5]{MR0434929}), where $\varphi(q)=q \prod_{p|q}\left(1-1/p\right)$ \marginpar{\textcolor{red}{Removed the multiplication dot ``$\cdot$'' between $q$ and $\prod_{p|q}$.}} denotes Euler totient function (see, e.g., \cite[\S 2.3 \& \S 2.5]{MR0434929}) \marginpar{\textcolor{red}{Added Section~\S 2.5.}}. In addition, if $\chi$ is primitive with $q>1$, the associated Dirichlet L-function satisfies the following functional equation, which holds for all complex numbers $s$ (see, e.g., \cite[\S 10, Corollary 10.9]{MR2378655}, \cite[Page 47]{MR1216135}, \cite{MR0285499}, \cite[\S 12.10]{MR0434929} and \cite[\S 5.4]{MR2376618}):
\begin{equation}\label{e:Fun_Equ_Dir_L-F}
L(1-s,\chi)=\frac{2\Gamma(s)\cos\left(\frac{\pi}{2}(s-\delta)\right)}{i^{\delta}(2\pi)^s}\frac{G(1,\chi)L(s,\overline{\chi})}{q^{1-s}}.
\end{equation}
Here $\delta=0$ if $\chi$ is even (i.e. $\chi(-1)=1$), $\delta=1$ if $\chi$ is odd (i.e. $\chi(-1)=-1$), and for $\alpha\in\mathbb{C}$,
\begin{equation}\label{e:Gauss_Sum_Fun}
G(\alpha,\chi)=\sum_{k=1}^{q}\chi(k)e^{\frac{2i\pi k\alpha}{q}}
\end{equation}
denotes the Gauss sum function associated with $\chi$ (see, e.g. \cite[Page 262]{MR0434929} and \cite[Page 378]{MR2860694}). This functional equation plays a crucial role in the analytic properties of Dirichlet L-functions, including their symmetry and behavior at integer values. The Dirichlet L-function and the Hurwitz zeta function are closely linked through the fundamental identity
\begin{equation}\label{e:LsChi_Hur}
L(s,\chi)=q^{-s}\sum_{a=1}^{q}\chi(a)\zeta(s,\frac{a}{q});
\end{equation} 
see, e.g., \cite[Page 249]{MR0434929}, \cite[Page 41]{MR1216135} and \cite[Page 71]{MR0606931}.

Euler's distinguished contribution to the field of mathematics includes the demonstration of an elegant formula in \S 178 of his \textit{Introductio in Analysin Infinitorum} \cite{MR1841793}, published in 1748. This formula, widely regarded as one of his most exceptional achievements, is the partial fraction expansion of the cotangent function: 
\begin{equation}\label{e:ParFraExpCot}
\pi\cot(\pi w)=\frac{1}{w}+\sum_{n=1}^{+\infty}\left(\frac{1}{w+n}+\frac{1}{w-n}\right),\ w\in\mathbb{C}\backslash\mathbb{Z}.
\end{equation}
In \cite[\S 23]{MR2569612}, Herglotz discovered an alternative method for proving this identity. His ingenious reasoning involved a remarkably simple technique that has since been referred to as the ``Herglotz trick''. He demonstrated that both functions on either side exhibit a comprehensive and common set of robust properties that justify the conclusion that they are identical. 

In \cite[Proposition 1]{MR1222534}, Louboutin using \eqref{e:LsChi_Hur} and \eqref{e:ParFraExpCot}, derived a closed expression for the value $ L(1,\chi) $ whenever $\chi$ is an odd non-principal Dirichlet character modulo $q\geq2$, not necessarily primitive:
\begin{equation}\label{e:L1Chi}
L(1,\chi)=\frac{\pi}{2q}\sum_{a=1}^{q-1}\chi(a)\cot\left(\frac{\pi a}{q}\right).
\end{equation}
In another paper \cite[Proposition.3.(1)]{MR1874365}, he generalizes \eqref{e:L1Chi} by providing a formula for the values of $L(k, \chi)$, where $\chi$ satisfies $\chi(-1) = (-1)^k$, i.e., $\chi$ and $k$ have the same parity, with $q \geq 3$ and $k \geq 2$:
\begin{equation}\label{e:LkChiCot}
L(k,\chi)=\frac{(-1)^{k-1}\pi^{k}}{2q^{k}(k-1)!}\sum_{a=1}^{q-1}\chi(a)\cot^{(k-1)}\left(\frac{\pi a}{q}\right).
\end{equation}
His proof involves applying \eqref{e:LsChi_Hur} and using the relation between the Hurwitz zeta function and the higher-order derivatives of the function $\cot(\pi w)$, $w\in\mathbb{C}\backslash\mathbb{Z}$. By the end of this article in Section~\ref{ss:Exp2}, an exact formula for these derivatives is derived, expressed in terms of $\cos(\pi w)$, $\sin(\pi w)$, and Eulerian numbers.
  
In \cite[Pages 336, 337]{MR2378655}, Exercises 14 and 15 pose questions about explicit formulas for the values $L(2p,\chi)$ and $L(2p+1,\chi) $, where $p$ is a positive integer and $\chi$ is a primitive character modulo $ q\geq2$:
\begin{equation}\label{e:EX14Form}
L(2p,\chi)=\frac{(-1)^{p+1}2^{2p-1}\pi^{2p}G(1,\chi)}{q(2p)!}\sum_{a=1}^{q}\overline{\chi}(a)B_{2p}\left(\frac{a}{q}\right),\text{ if }\chi(-1)=1,
\end{equation}
\begin{equation}\label{e:EX15Form}
L(2p+1,\chi)=\frac{i(-1)^{p}2^{2p}\pi^{2p+1}G(1,\chi)}{q(2p+1)!}\sum_{a=1}^{q}\overline{\chi}(a)B_{2p+1}\left(\frac{a}{q}\right),\text{ if }\chi(-1) = -1.
\end{equation}
Here $B_n(x)$ is the $n$-th Bernoulli polynomial, and $\overline{\chi}$ represents the complex conjugate of the character $\chi$. The specific case of a closed formula for the value $ L(1,\chi)$ is detailed in \cite[Theorem 9.9]{MR2378655}.

As a consequence of \cite[Theorem 1]{MR2860694} \marginpar{\textcolor{red}{Replaced ``Theorem 1 in \cite[Page 380]{MR2860694}'' with ``\cite[Theorem 1]{MR2860694}''.}}, Alkan established a new exact formula for $L(k, \chi)$, where $k\in\mathbb{Z}_{\geq 1}$ and $\chi$ is a Dirichlet character modulo $q\geq2$ satisfying $\chi(-1)=(-1)^k$ (see, \cite[Page 380]{MR2860694}): \marginpar{\textcolor{red}{Added ``(see, [2, Page 380])''.}}
\begin{equation}\label{e:alkan_formula}
L(k,\chi)=\frac{(-1)^{k+1} i^{k} 2^{k-1}\pi^{k}}{q k! }\sum_{l=0}^{2\left\lfloor\frac{k}{2}\right\rfloor} \binom{k}{l}B_l\sum_{a=1}^{q}\left(\frac{a}{q}\right)^{k-l}G(a,\chi),
\end{equation}
where $B_n$ is the $n$-th Bernoulli number, and $\lfloor\cdot\rfloor$ denotes the floor function. Subsequently, under the assumption that $\chi$ is primitive, he derived the formula
\begin{equation}\label{eq:Alk_for}
L(k,\chi)=\frac{(-1)^{k+1}i^{k} 2^{k-1}\pi^{k}G(1,\chi)}{qk!}\sum_{a=1}^{q}\overline{\chi}(a)B_k\left(\frac{a}{q}\right),
\end{equation}
which serves as a unification of the aforementioned formulas \eqref{e:EX14Form} and \eqref{e:EX15Form}. Finally, by employing the generalized Bernoulli
numbers (see, e.g., \cite[Page 43]{MR1216135} and \cite[Page 441]{MR1697859}) given by
\begin{equation}\label{eq:gen_Ber_num}
B_{n,\chi}=q^{n-1}\sum_{a=1}^{q}\chi(a)B_n\left(\frac{a}{q}\right),
\end{equation}
Alkan elegantly recovers the classical formula (see, e.g., \cite[Page 443, Corollary 2.10]{MR1697859})
\begin{equation}\label{eq:L_k_Chi_cla_for}
L(k,\chi)=\frac{(-1)^{k+1} i^{k} 2^{k-1}\pi^{k}G(1,\chi)}{q^{k}k!}B_{k,\overline{\chi}}.
\end{equation}
This formulation encapsulates the relationship between Dirichlet L-functions and generalized Bernoulli numbers, providing a concise expression for $L(k,\chi)$ in terms of these number-theoretic objects.

In \cite{Lamgouni2024}, we conducted a detailed investigation of LC-functions, generalizing the Hurwitz zeta function and introducing an extended Hurwitz formula. This paper further explores this LC-function formula and reveals behaviors of certain LC-functions analogous to those of Dirichlet L-functions.

More precisely, throughout this paper, let $f$ be a real-analytic function at $0$, defined by the series
\begin{equation}\label{e:f_Ser}
f(t)=\sum_{n=0}^{+\infty}\frac{C_{f,n}}{n!}t^n,
\end{equation}
such that the series $\sum_{n=0}^{+\infty}P_{f,n}t^{n}$ has a non-zero radius of convergence $\rho_f$. Here, $P_{f,n}$ are the \textit{P-numbers} associated to $f$ defined through the exponential generating function
\begin{equation}\label{e:pf_Ser}
p_f(t):=\frac{e^{t}-1}{t}f(t)=\sum_{n=0}^{+\infty}\frac{P_{f,n}}{n!}t^n,
\end{equation}
and $C_{f,n}$ are the \textit{C-numbers} associated to $f$; see \cite[\S 2]{Lamgouni2024} for more details on C-numbers and P-numbers. As usual, denote by $\sigma$ the real part of the complex variable $s$. The LC-function associated to $f$ (see \cite[\S 5]{Lamgouni2024}) is defined for $\sigma>1$ by the Dirichlet series
\begin{equation}\label{e:LC_Fun_Ser}
L(s,f):=\sum_{n=n_f}^{+\infty} n^{(-s,f)},
\end{equation}
where $n_f:=\lfloor 1/\rho_f\rfloor+1$, and the bivariate complex function $z^{(s,f)}$ (see \cite[\S 4]{Lamgouni2024}) defined for all $s\in\mathbb{C}$ and all $z\in\mathbb{C}\setminus\mathbb{R}_{\leq 0}$ with $\abs{z}>1/\rho_f$ by
\begin{equation}\label{e:Pow_szf1}
z^{(s,f)}:=z^{s}\sum_{k=0}^{+\infty}\binom{s}{k}P_{f,k}\left(\frac{1}{z}\right)^{k},
\end{equation}
is a generalization of the complex exponentiation $z^{s}$.

The LC-function $L(s,f)$ extends to an entire function when $p_f(0) = 0$, otherwise to a meromorphic function with a simple pole at $s=1$ with residue $p_f(0)$ (see, \cite[\S 5.4]{Lamgouni2024}). \marginpar{\textcolor{red}{Removed ``e.g.'' before ``see''.}} Moreover, for $\sigma<0$, this function satisfies the following equation, henceforth referred to as the \textit{LC-formula} (see \cite[\S 5.6]{Lamgouni2024}):
\begin{equation}\label{e:LC_For}
L(1-s,f)=\frac{\Gamma(s)}{(2\pi)^{s}}\left(e^{-\frac{i\pi s}{2}}F(s,f_{(-2i\pi)})+e^{\frac{i\pi s}{2}}F(s,f_{(2i\pi)})\right).
\end{equation}
Here $f_{(\alpha)}$, with $\alpha\in\mathbb{C}$, denotes the function 
\begin{equation}\label{e:f_Alp}
f_{(\alpha)}(t):=\frac{e^{\alpha t}-1}{\alpha(e^{t}-1)}f(\alpha t),
\end{equation}
and $F(s,f)$ is the \textit{FC-function} associated to $f$ defined by
\begin{equation}\label{e:FC_Fun_Con_Int}
F(s,f):=\frac{\Gamma(1-s)}{2i\pi}\int_{\mu_f}\frac{z^{(s-1,f)}}{e^{-z}-1}\, dz,
\end{equation}
where $\mu_f$ is the Hankel path used in the contour integration. The FC-function and its associated Hankel contour are described in detail in Section~\ref{ss:FC_Fun}.

The function $p_f$, stated in \eqref{e:pf_Ser} and henceforth referred to as the \textit{modulator} of the LC-function $L(s, f)$, is essential in this investigation. Later in this paper, we show that characteristics of $p_f$ significantly influence the properties of $L(s, f)$, just as characteristics of a Dirichlet character $\chi$ influence the properties of the L-function $L(s,\chi)$.

Choosing $a=1$ in the Hurwitz formula \eqref{e:Hur_For} and recognizing that $\zeta(s,1)=\zeta(s)$ and $F(1,s)=F(-1,s)=\zeta(s)$, Euler's formula $e^{i\pi s/2}+e^{-i\pi s/2}=2\cos(\pi s/2)$ can be applied to simplify \eqref{e:Hur_For} to the famous functional equation of the Riemann zeta function (see, e.g., the proof of Theorem 12.7 in \cite{MR0434929}):
\begin{equation}\label{e:Fun_Equ_Rie_Zet}
\zeta(1-s)=2(2\pi)^{-s}\Gamma(s)\cos\left(\frac{\pi s}{2}\right)\zeta(s).
\end{equation}
Noting the striking analogy between the LC-formula \eqref{e:LC_For} and Hurwitz formula \eqref{e:Hur_For}, a natural question arises: Can the LC-formula be simplified in a similar manner, and which class of functions $ f $ permits this simplification? This question is crucial because such a simplification will provide deeper insights into the analytic properties of a subclass of LC-functions, which, as a generalization, exhibit behavior similar to the Riemann zeta function. We affirmatively answer this question by focusing on cases where the modulator $ p_f $ exhibits even or odd symmetry. The main result of this paper, presented in Theorem~\ref{t:Mai_The_v1} below, demonstrates that the LC-formula for an LC-function with either an even or odd modulator simplifies to a functional equation remarkably similar to that of a Dirichlet L-function \eqref{e:Fun_Equ_Dir_L-F} with a primitive character.

\begin{theorem}[cf. Theorem~\ref{t:Mai_The_v2}]
\label{t:Mai_The_v1}
Let $L(s,f)$ be an LC-function with modulator $p_f$ of even or odd parity. For all $s\in\mathbb{C}\setminus\{0\}$,
\begin{equation}\label{e:Fun_Equ_LC_Fun}
L(1-s,f)=\frac{2i^{\kappa}\Gamma(s)\cos\left(\frac{\pi}{2}(s-\kappa)\right)}{(2\pi)^{s}}F(s,f_{(2i\pi)}).
\end{equation}
Here $\kappa=0$ if $p_f$  is even and $\kappa=1 $ if $p_f$ is odd. 
\end{theorem}

The paper is organized as follows. In Section~\ref{s:LC_Fun_LC_For}, we revisit LC-functions and their LC-formula in greater detail. In Section~\ref{s:Fun_Equ_LC_Fun_EOM}, we prove the main result, Theorem~\ref{t:Mai_The_v2}, and subsequently, in Section~\ref{s:Spe_Val_LC_fun_FC_fun}, we use this result to establish two general formulas presented in Theorem~\ref{t:Spe_Val_LC_fun_FC_fun}: The first one concerns the values at positive even integers of LC-functions with an even modulator, and the second addresses the values at positive odd integers of LC-functions with an odd modulator. In Section~\ref{s:Val_Ser_rep_LC_Fun_VM}, we demonstrate that the series representation given in \eqref{e:LC_Fun_Ser} for LC-functions is valid in the half-plane $\sigma>0$, provided that its modulator $p_f(t)$ vanishes at $0$. This bears a notable analogy to the series representation of Dirichlet L-functions with non-principal character, which are also defined within the half-plane $\sigma>0$ as the Gauss sum function \eqref{e:Gauss_Sum_Fun} vanishes at $0$ (i.e., $G(0,\chi)=\sum_{k=1}^{q}\chi(k)=0$). Finally, in Section~\ref{s:App}, we present two illustrative examples to showcase the applicability of the tool established in Section~\ref{s:Spe_Val_LC_fun_FC_fun}. The first example provides an explicit formula for the values at positive even integers of the LC-function attached to the even modulator $ t\mapsto\cos(wt) $, and the second example provides an explicit formula for the values at positive odd integers of the LC-function attached to the odd modulator $t\mapsto\sinh(wt)$. Here, $w$ is a real parameter satisfying $\abs{w}<1$. Concurrently, we provide general closed-form formulas for the two integrals in terms of elementary functions:
\[
\int_0^{+\infty}\frac{t^{2p-1}\cos(wt)}{e^t-1}\,dt=\frac{(-1)^p}{2}\left(\frac{(2\pi)^{2p}e^{2\pi w}A_{2p-1}(e^{2\pi w})}{(1-e^{2\pi w})^{2p}}-\frac{(2p-1)!}{w^{2p}}\right),
\]
and
\[
\int_0^{+\infty}\frac{t^{2p}\sinh(wt)}{e^t-1}\,dt=\frac{1}{2}\left(\frac{(2i\pi)^{2p+1}e^{2i\pi w}A_{2p}(e^{2i\pi w})}{(1-e^{2i\pi w})^{2p+1}}+\frac{(2p)!}{w^{2p+1}}\right).
\]
Here, $p$ is a positive integer and $A_n(x)$ denotes the Eulerian polynomials. Furthermore, we derive a general closed-form formula, in terms of $\cos(\pi w)$, $\sin(\pi w)$, and the Eulerian numbers for the functions whose partial fraction expansion takes the form
\[
\frac{1}{w^{p}}+\sum_{n=1}^{+\infty}\left(\frac{1}{(w-n)^{p}}+\frac{1}{(w+n)^{p}}\right).
\]
For details on Eulerian numbers and polynomials, we refer readers to, e.g., \cite{MR2744266} and \cite[\S 1.3 \& \S 1.4]{MR3408615}. Through analytic continuation, the results above can be extended to complex $w$, as specified in Section~\ref{s:App}.

As a noteworthy consequence of our main result, we find that Euler's formula \eqref{e:ParFraExpCot} naturally follows from the value at $1$ of the LC-function associated to the function $\phi(t)=t\sinh(wt)/(e^t-1)$. Specifically, we state that:
\[
L(1,\phi)=-\frac{1}{2w}+\frac{\pi}{2}\cot(\pi w).
\]
This expression, to be derived later in Section~\ref{ss:Exp2}, can be compared to the formula \eqref{e:L1Chi} for $L(1,\chi)$. This analogy, along with others observed throughout the article, suggests a profound potential connection between the L-functions and LC-functions associated with modulators of even or odd parity. The precise nature of the parallel between these two classes of functions, however, remains an open question.

\section{Background material}\label{s:LC_Fun_LC_For}
To ensure that our discussion is self-contained and to establish the foundation for the main contributions of this paper, this section reviews the necessary definitions and results related to LC-functions, FC-functions, and the LC-formula. For more details, we refer readers to \cite{Lamgouni2024}.

\subsection{Generalized complex exponentiation}\label{ss:Gen_Com_Exp}
\begin{figure}
\includegraphics{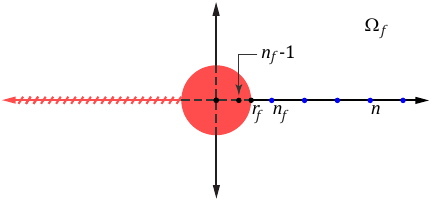}
\caption{Let $s\in \mathbb{C}$ be fixed. The generalized complex exponentiation function $z\mapsto z^{(s,f)}$  is well-defined on $\Omega_f$, particularly for all integers $n\geq n_f$.}
\label{f:Omega_f}
\end{figure}

As defined previously in \eqref{e:Pow_szf1}, the bivariate complex function $z^{(s,f)}$ is given for all $s\in\mathbb{C}$ and $z\in\Omega_f$ (see Figure~\ref{f:Omega_f}) by
\begin{equation}\label{e:Pow_szf2}
z^{(s,f)}=\sum_{k=0}^{+\infty}\binom{s}{k}P_{f,k}z^{s-k},
\end{equation}
where
\begin{equation}\label{e:Omega}
\Omega_f:=\left\{\alpha\in\mathbb{C}\mid\alpha\notin\mathbb{R}_{\leq 0}\land\abs{\alpha}>r_f\right\}.
\end{equation}
Here, $r_f:=1/\rho_f$, with $r_f=0 $ if $\rho_f=\infty$. This function generalizes the complex exponentiation $z^{s}=e^{s\log z}$, where $\log$ denotes the principal branch of the logarithm (used consistently throughout this article). For more details about the function $z^{(s,f)}$, see \cite[\S 4]{Lamgouni2024}. Within the domain $\Omega_f$, the smallest positive integer, denoted $n_f$, is defined by the equation
\begin{equation}\label{e:nf}
n_f:=\lfloor r_f\rfloor+1.
\end{equation}
This ensures that the function of $s$, $n^{(s,f)}$, for $n\in\mathbb{Z}_{\geq0}$, is well-defined exclusively for all integers $n\geq n_f$ (see Figure~\ref{f:Omega_f}), with $s$ being any complex number.
Building upon the function $(z,s)\mapsto z^{(s,f)}$, we defined two fundamental functions that extend its applications: LC-functions and FC-functions.

\subsection{LC-functions}\label{ss:LC_Fun}
This subsection presents various representations of the LC-function, each offering unique insights into its properties. The LC-function associated to $f$ has a series representation presented in \eqref{e:LC_Fun_Ser} that converges absolutely in the half-plane $\sigma>1$. This function can alternatively be expressed using the following integral representation, which is valid for $\sigma>1 $ (see \cite[Theorem 5.2]{Lamgouni2024}),
\begin{equation}\label{e:LC_Fun_Int}
L(s,f)=\frac{1}{\Gamma(s)}\int_0^{+\infty}\frac{t^{s-1}e^{(1-n_f)t}p_f(-t)}{e^t-1}\,dt=\frac{1}{\Gamma(s)}\int_0^{+\infty} t^{s-2}e^{-n_ft}f(-t)\,dt.
\end{equation}
To provide further insight, the integral representation above can be converted into the contour integral representation (see \cite[Theorem 5.3]{Lamgouni2024})
\begin{equation}\label{e:LC_Fun_Con_Int}
L(s,f)=\frac{\Gamma(1-s)}{2i\pi}\int_{\mu}\frac{z^{s-1}e^{n_fz}p_f(z)}{1-e^{z}}\,dz,
\end{equation}
where $\mu$ is a Hankel contour counterclockwise around the negative real axis, and the principal branch is used for the complex exponentiation $z^{s-1}$. It should be noted that the modulator $p_f$ is an entire complex function; see \cite[Theorem 4.1.(ii)]{Lamgouni2024}. 
\begin{remark}\label{r:LC_Fun_Pol}
According to \cite[Theorem 5.4]{Lamgouni2024}, if $p_f(0)\neq0$ (i.e., $f(0)\neq0$), the contour integral representation provides an analytic continuation of $L(s, f)$ to the whole complex plane, excluding the point $1$, where it presents a simple pole with a residue equal to $p_f(0)$. On the other hand, when $p_f(0)=0$, $L(s, f)$ becomes an entire function. This behavior reveals the first notable analogy between Dirichlet L-functions and LC-functions. Specifically, the Dirichlet L-function $L(s,\chi_0)$ for the principal character $\chi_0$ modulo $q$ (i.e., $G(0,\chi_0)=\sum_{k=1}^{q}\chi_0(k)\neq0$) has a simple pole at $s=1$ with residue $G(0,\chi_0)/q$ \marginpar{\textcolor{red}{Corrected: replaced $\chi$ with $\chi_0$.}}. In contrast, the L-function $L(s,\chi)$ for a non-principal character $\chi$ (i.e., $G(0,\chi)=0$) is an entire function.
\end{remark}

\subsection{FC-functions}\label{ss:FC_Fun}

\begin{figure}
\includegraphics{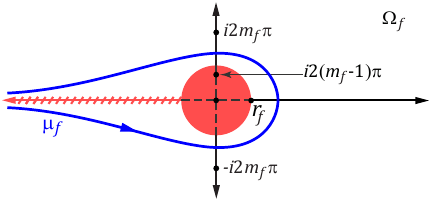}
\caption{The Hankel contour $\mu_f$ is oriented counterclockwise around the negative real axis and the closed disk $\abs{z}\leq r_f$, so as not to encircle any points of discontinuity of the function $1/(e^{-z} - 1)$ that lie within the region $\Omega_f$. The integer $m_f=\lfloor r_f/(2\pi)\rfloor+1$ is defined as the smallest positive integer $m$ such that $i2m\pi\in\Omega_f$.}
\label{f:Han_Con_FC_Fun}
\end{figure}

We now introduce the FC-function, another fundamental concept in this study, defined through a contour integral representation \eqref{e:FC_Fun_Con_Int}. Here, the Hankel path $ \mu_f $ is oriented counterclockwise around the negative real axis and the closed disk $ \abs{z} \leq r_f $. It does not encircle any points of discontinuity of the function $ 1/(e^{-z} - 1) $ that lie within the region $ \Omega_f $, as shown in Figure~\ref{f:Han_Con_FC_Fun}. The value $ m_f $ is the smallest positive integer $ m $ satisfying $ i2m\pi \in \Omega_f $. It is given by
\begin{equation}\label{e:mf}
m_f=\left\lfloor\frac{r_f}{2\pi}\right\rfloor+1.
\end{equation}
The function $F(s,f)$ is confirmed in \cite[Theorem 5.6]{Lamgouni2024} to be analytic for $s\neq 1,\ 2,\ 3,\dotsc$. However, this analyticity does not necessarily imply that this function has poles at these specific points; the behavior of $F(s,f)$ at these integers remains to be further analyzed, depending on $f$. The integral in \eqref{e:FC_Fun_Con_Int} is itself an entire function of $s$.

\subsection{LC-formula}\label{ss:LC_For}
As introduced in \eqref{e:LC_For}, the LC-formula unifies the two functions LC and FC through a single equation that holds for $s\neq0$ \marginpar{\textcolor{red}{Corrected: replaced $\sigma\neq0$ with $s\neq0$.}}. Indeed, according to Remark~\ref{r:LC_Fun_Pol}, the LC-function $L(s,f)$ is known to be analytic throughout the entire complex plane except, possibly, at the point $1$---based on whether $p_f(0)$ vanishes or not. This fact establishes that the singularities $\dotsc$, $-3$, $-2$, $-1$, $1$, $2$, $3$, $\dotsc$ of the right-hand side of \eqref{e:LC_For} are removable.

\begin{figure}
\includegraphics{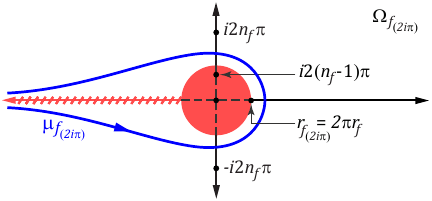}
\caption{This figure illustrates the Hankel contour $\mu_{f_{(2i\pi)}}$, shown to be equivalent to $\mu_{f_{(-2i\pi)}}$.}
\label{f:Han_Con_FC_Fun_2iPi}
\end{figure}

Now, we add a brief clarification on the Hankel contour related to the FC-function $ F(s,f_{(\omega)}) $, where $\omega=\pm 2i\pi$. From Equations \eqref{e:pf_Ser} and \eqref{e:f_Alp}, we deduce that $p_{f_{(\omega)}}(t)=p_f(\omega t)$. Consequently,
\[
\sum_{n=0}^{+\infty}\frac{P_{f_{(\omega)},n}}{n!}t^n=\sum_{n=0}^{+\infty}\frac{\omega^{n}P_{f,n}}{n!}t^n.
\]
Thus, for all $n\geq0$,
\begin{equation}\label{e:Pfwn_wnPfn}
P_{f_{(\omega)},n}=\omega^{n}P_{f,n}.
\end{equation}
It follows then that the radius of convergence of the series $\sum_{n=0}^{+\infty}P_{f_{(\omega)},n}t^{n}$ is $ \rho_{f_{(\omega)}}=\rho_f / \abs{\omega}$. Equivalently, 
\begin{equation}\label{e:rfw_2Pirf}
r_{f_{(\omega)}}=1/\rho_{f_{(\omega)}}=2\pi r_f.
\end{equation}
Thus, $r_{f_{(-2i\pi)}}=r_{f_{(2i\pi)}}$, and by Equations \eqref{e:mf} and \eqref{e:nf},
\[
m_{f_{(-\omega)}}=m_{f_{(\omega)}}=\lfloor r_f\rfloor+1=n_f.
\]
Hence, the Hankel paths $\mu_{f_{(-2i\pi)}}$ and $\mu_{f_{(2i\pi)}}$ are equivalent; see Figure~\ref{f:Han_Con_FC_Fun_2iPi}. This equivalence simplifies the analysis of the FC-functions $F(s, f_{(\omega)})$, allowing for a unified treatment in subsequent calculations within the proof of the Main Theorem.

To conclude this section, we provide a formula for the values of the LC-functions at non-positive integers, which is later utilized to state some identities in Theorem~\ref{t:Spe_Val_LC_fun_FC_fun}. In \cite[Theorem 5.5]{Lamgouni2024}, it is established that for all integers $p\geq0$,
\begin{equation}\label{e:L_min_p}
L(-p,f)=-\frac{C_{f,p+1}(n_f)}{p+1}=(-1)^p\frac{C_{\underline{f},p+1}(1-n_f)}{p+1}.
\end{equation}
Here $\underline{f}(t)$ denotes the function defined by (see \eqref{e:f_Alp})
\begin{equation} \label{e:f_undlin}
\underline{f}(t):=f_{(-1)}=e^{-t}f(-t),
\end{equation}
and $C_{f,n}(x)$ are the \textit{C-polynomials} associated to $f$ (introduced in \cite[\S 2]{Lamgouni2024}), defined by the exponential generating function
\begin{equation}\label{e:C_Pol_Ser}
f(t)e^{xt}=\sum_{n=0}^{+\infty}\frac{C_{f,n}(x)}{n!}t^n.
\end{equation}

\section{The Main Theorem}\label{s:Fun_Equ_LC_Fun_EOM}

\subsection{Functional equation for LC-functions with even or odd modulator}\label{ss:Fun_Equ_LC_Fun_EOM}
In this subsection, we present a reformulation of the Main Theorem, Theorem~\ref{t:Mai_The_v1}, by dividing the functional equation \eqref{e:Fun_Equ_LC_Fun} into two cases based on the parity of the modulator of the LC-function. We then proceed to provide the proof.

\begin{theorem}\label{t:Mai_The_v2}
\begin{enumerate}[label=(\roman*)]
\item Let $L(s,f)$ be an LC-function with even modulator $p_f$. For all $ s\in\mathbb{C}\setminus\{0\} $,
\begin{equation}\label{e:Fun_Equ_LC_Fun_Cos}
L(1-s,f)= 2(2\pi)^{-s}\Gamma(s)\cos\left(\frac{\pi s}{2}\right)F(s,f_{(2i\pi)}).
\end{equation}
Using the definition \eqref{e:FC_Fun_Con_Int} of an FC-function and applying Euler's reflection formula, we obtain the following form of \eqref{e:Fun_Equ_LC_Fun_Cos}
\begin{equation}\label{e:Fun_Equ_LC_Fun_SInt}
L(1-s,f)=\frac{(2\pi)^{-s}}{2i\sin\left(\frac{\pi s}{2}\right)}\int_{\mu_{f_{(2i\pi)}}}\frac{z^{(s-1,f_{(2i\pi)})}}{e^{-z}-1}\,dz.
\end{equation}
\item Let $L(s,f)$ be an LC-function with odd modulator $p_f$. For all $s\in\mathbb{C}$,
\begin{equation}\label{e:Fun_Equ_LC_Fun_Sin}
L(1-s,f)=2i(2\pi)^{-s}\Gamma(s)\sin\left(\frac{\pi s}{2}\right)F(s,f_{(2i\pi)}).
\end{equation}
Similarly, we can rewrite \eqref{e:Fun_Equ_LC_Fun_Sin} as
\begin{equation}\label{e:Fun_Equ_LC_Fun_CInt}
L(1-s,f)=\frac{(2\pi)^{-s}}{2\cos\left(\frac{\pi s}{2}\right)}\int_{\mu_{f_{(2i\pi)}}}\frac{z^{(s-1,f_{(2i\pi)})}}{e^{-z}-1}\,dz.
\end{equation}
\end{enumerate}
\end{theorem}

\begin{remark}
\label{r:Ana_Con_LC_Fun}
(i) Assume that the modulator $p_f$ is even. If $p_f(0)\neq0$, then according to Remark~\ref{r:LC_Fun_Pol}, the LC-function $L(s,f)$ is holomorphic on the entire complex plane, except at $s=1$, where it has a simple pole with a residue equal to $p_f(0)$. Therefore, the singularities on the right-hand side of \eqref{e:Fun_Equ_LC_Fun_Cos} and \eqref{e:Fun_Equ_LC_Fun_SInt} are removable, except at $s=0$. However, if $p_f(0)=0$, then by the same Remark, $L(s, f)$ is an entire function. Thus, all singularities on the right-hand side of these two equations are removable.

(ii) If the modulator $p_f$ is odd, then $p_f(0)=0$. Therefore, according to Remark~\ref{r:LC_Fun_Pol}, $L(s,f)$ is entire. Consequently, all singularities on the right-hand side of \eqref{e:Fun_Equ_LC_Fun_Sin} and \eqref{e:Fun_Equ_LC_Fun_CInt} are removable.
\end{remark}

\begin{proof}[Proof of the Main Theroem]
Let $L(s,f)$ be an LC-function associated to a given function $f$, whose modulator $p_f$ is either an even or odd function. Initially, we establish that if $p_f$ is even, the corresponding FC-functions $F(s,f_{(-2i\pi)})$ and $F(s,f_{(2i\pi)})$ are identical. Conversely, if $p_f$ is odd, these two FC-functions are opposite.

Based on \eqref{e:Pow_szf2} and \eqref{e:Pfwn_wnPfn}, one can write
\begin{equation}\label{e:Pow_szfw}
z^{(s,f_{(-2i\pi)})}=\sum_{n=0}^{+\infty}\binom{s}{n}(-2i\pi)^nP_{f,n}z^{s-n}.
\end{equation}
Now, assuming that $p_f$ is even, all the coefficients $P_{f,2n+1}$ in \eqref{e:pf_Ser} are zero. Hence, we obtain
\[
z^{(s,f_{(-2i\pi)})}=\sum_{n=0}^{+\infty}\binom{s}{2n}(2i\pi)^{2n}P_{f,2n}z^{s-2n}=z^{(s,f_{(2i\pi)})}.
\]
Therefore, by \eqref{e:FC_Fun_Con_Int}, since the Hankel paths $\mu_{f_{(-2i\pi)}}$ and $\mu_{f_{(2i\pi)}}$ are equivalent (see Figure~\ref{f:Han_Con_FC_Fun_2iPi}),
\[
F(s,f_{(-2i\pi)})=\frac{\Gamma(1-s)}{2i\pi}\int_{\mu_{f_{(-2i\pi)}}}\frac{z^{(s-1,f_{(2i\pi)})}}{e^{-z}-1}\,dz=F(s,f_{(2i\pi)}).
\]

On the other hand, if we assume $p_f$ is odd, then all the coefficients $P_{f,2n}$ vanish. Consequently, \eqref{e:Pow_szfw} reduces to
\[
z^{(s,f_{(-2i\pi)})}=-\sum_{n=0}^{+\infty}\binom{s}{2n+1}(2i\pi)^{2n+1} P_{f,2n+1}z^{s-2n-1}=-z^{(s,f_{(2i\pi)})}.
\]
Thus, similarly,
\[
F(s,f_{(-2i\pi)})=-\frac{\Gamma(1-s)}{2i\pi}\int_{\mu_{f_{(-2i\pi)}}}\frac{z^{(s-1,f_{(2i\pi)})}}{e^{-z}-1}\,dz= -F(s,f_{(2i\pi)}).
\]
In order to unify the proofs of formulas \eqref{e:Fun_Equ_LC_Fun_Cos} and \eqref{e:Fun_Equ_LC_Fun_Sin}, we introduce the parameter
\begin{align*}
\kappa=\kappa(f)=\begin{cases} 
0&\text{if } p_f \text{ is even}, \\
1&\text{if } p_f \text{ is odd}.
\end{cases}
\end{align*}
With this notation, the identity $ F(s,f_{(-2i\pi)})=(-1)^{\kappa}F(s,f_{(2i\pi)}) $ holds regardless of whether $ p_f $ is even or odd. Thus, the LC-formula \eqref{e:LC_For} becomes
\[
L(1-s,f)=\frac{\Gamma(s)}{(2\pi)^{s}}\left((-1)^{\kappa}e^{\frac{-i\pi s}{2}}+e^{\frac{i\pi s}{2}}\right)F(s,f_{(2i\pi)}).
\]
Utilizing the identity $(-1)^{\kappa}=i^{2\kappa}$ and substituting it into the formula, we obtain
\[
L(1-s,f)=\frac{i^{\kappa}\Gamma(s)}{(2\pi)^{s}}\left(i^{\kappa}e^{\frac{-i\pi s}{2}}+i^{-\kappa}e^{\frac{i\pi s}{2}}\right)F(s,f_{(2i\pi)}).
\]
Subsequently, we write
\[
i^{\kappa}e^{\frac{-i\pi s}{2}}+i^{-\kappa}e^{\frac{i\pi s}{2}}=e^{\frac{i\pi}{2}(s-\kappa)}+e^{-\frac{i\pi}{2}(s- \kappa)}.
\]
Thus, we arrive at the final simplified expression
\[
L(1-s,f)=\frac{2i^{\kappa}\Gamma(s)}{(2\pi)^{s}}\cos\left(\frac{\pi}{2}(s-\kappa)\right)F(s,f_{(2i\pi)}).
\]
\end{proof}

\subsection{Riemann zeta function as an LC-function with even modulator}\label{ss:Riem_Zeta_LC_func}
The Riemann zeta function is the LC-function, $L(s,\beta)$, associated to the function $\beta(t):=t/(e^{t}-1)$. Indeed, from \eqref{e:pf_Ser}, we have 
\[
p_{\beta}(t)=1=\sum_{n=0}^{+\infty}\frac{P_{\beta,n}}{n!}t^n,
\]
which implies that all the P-numbers $P_{\beta,n}$ are zero except for $P_{\beta,0}$, which equals $1$. Consequently, the series $\sum_{n=0}^{+\infty}P_{\beta,n} t^{n}=1$ has an infinite radius of convergence, i.e., $\rho_{\beta}=\infty$. From \eqref{e:nf}, we have $n_{\beta} = 1$. Furthermore, according to \eqref{e:Pow_szf2}, for all $z\in\Omega_{\beta}=\mathbb{C}\setminus\mathbb{R}_{\leq 0}$ and $s\in\mathbb{C}$, we have $z^{(s,\beta)}=z^s$. Since $p_{\beta}$ is the constant function $1$, $L(s,\beta)$ is an LC-function with even modulator. It is defined for $\sigma>1$ according to \eqref{e:LC_Fun_Ser} by the series
\[
L(s,\beta)=\sum_{n=n_{\beta}}^{+\infty}n^{(-s,\beta)}=\zeta(s).
\]
We note that, according to \eqref{e:f_Alp}, $\beta_{(\alpha)}=\beta$ for all $\alpha\in\mathbb{C}^{*}$. Therefore, from \eqref{e:FC_Fun_Con_Int}, it follows that the FC-function $F(s,\beta_{(2i\pi)})$ also coincides with the Riemann zeta function, namely,
\[
F(s,\beta_{(2i\pi)})=F(s,\beta)=\frac{\Gamma(1-s)}{2i\pi}\int_{\mu_{\beta}}\frac{z^{s-1}}{e^{-z}-1}\,dz=\zeta(s);
\]
see, e.g., \cite[Theorem 12.3]{MR0434929} with $a=1$. As a result, the functional equation \eqref{e:Fun_Equ_LC_Fun_Cos} for the LC-function $L(s,\beta)$ matches the functional equation for the Riemann zeta function \eqref{e:Fun_Equ_Rie_Zet}.

\section{Specific values of LC-functions and FC-functions with even or odd modulators}\label{s:Spe_Val_LC_fun_FC_fun}

The functional equation for an LC-function $L(s,f)$, as established in Theorem~\ref{t:Mai_The_v2}, yields notable formulas depending on the parity of the modulator $p_f$. When $p_f$ is even, these formulas are employed to evaluate $L(2p,f)$ and $F(2p,f_{(2i\pi)})$. Conversely, when $p_f$ is odd, they are used to evaluate $L(2p+1,f)$ and $F(2p+1,f_{(2i\pi)})$. In both cases, $p$ represents a non-negative integer.

\begin{theorem} \label{t:Spe_Val_LC_fun_FC_fun}
\begin{enumerate}[label=(\roman*)]
\item If $p_f$ is even, we have for every integer $p\geq0$
\begin{equation}\label{e:L2p_Eve}
L(2p,f)=\frac{(-1)^{p}(2\pi)^{2p}}{4i\pi}\int_{\mu_{f_{(2i\pi)}}}\frac{z^{(-2p,f_{(2i\pi)})}}{e^{-z} - 1}\,dz,
\end{equation}
\begin{equation}\label{e:F2p_Eve}
F(2p,f_{(2i\pi)})=\frac{(-1)^{p+1}C_{\underline{f},2p}(1-n_f)}{2(2p)!}(2\pi)^{2p}.
\end{equation} 
\item If $p_f$ is odd, we have for every integer $p\geq0$
\begin{equation}\label{e:L2p1_Od}
L(2p+1,f)=\frac{(-1)^p(2\pi)^{2p+1}}{4\pi}\int_{\mu_{f_{(2i\pi)}}}\frac{z^{(-2p-1,f_{(2i\pi)})}}{e^{-z}-1}\,dz,
\end{equation}
\begin{equation}\label{e:F2p1_Od}
F(2p+1,f_{(2i\pi)})=\frac{(-1)^{p}C_{\underline{f},2p+1}(1-n_f)}{2i(2p+1)!}(2\pi)^{2p+1}.
\end{equation}
\end{enumerate}
\end{theorem}

\begin{proof}
(i) Equation \eqref{e:L2p_Eve} follows by substituting $1-2p$ for $s$ in \eqref{e:Fun_Equ_LC_Fun_SInt}. For $p>0$, \eqref{e:F2p_Eve} is derived from \eqref{e:Fun_Equ_LC_Fun_Cos} by first replacing $s$ with $2p$, and then applying \eqref{e:L_min_p}. For the particular case $p=0$, we combine \eqref{e:Fun_Equ_LC_Fun_Cos} and \eqref{e:LC_Fun_Con_Int} to obtain the formula
\[
\frac{1}{2i\pi}\int_{\mu}\frac{z^{-s}e^{n_fz}p_f(z)}{1-e^{z}}\,dz=2(2\pi)^{-s}\cos\left(\frac{\pi s}{2}\right)F(s, f_{(2i\pi)}).
\]
According to \eqref{e:pf_Ser}, $p_f(z)/(1-e^{z})=-f(z)/z$. Hence,
\[
F(0,f_{(2i\pi)})=\frac{-1}{4i\pi}\int_{\mu}z^{-1}e^{n_fz}f(z)\,dz.
\]
We now apply Cauchy's residue theorem, noting from \eqref{e:C_Pol_Ser} that
\[
z^{-1}e^{n_fz}f(z)=\sum_{n=0}^{+\infty}\frac{C_{f,n}(n_f)z^{n-1}}{n!}.
\]
The contour $\mu$ can be reduced to a circle of sufficiently small radius centered at the origin. Thus, we obtain
\[
F(0,f_{(2i\pi)})=-\frac{C_{f,0}(n_f)}{2}=-\frac{C_{f,0}}{2}.
\]
Consequently, \eqref{e:F2p_Eve} is valid even for $p=0$, as $C_{\underline{f},0}(x)=C_{f,0}$, according to \eqref{e:C_Pol_Ser}, \eqref{e:f_undlin} and \eqref{e:f_Ser}.

(ii) Similar to the proof in (i), \eqref{e:L2p1_Od} follows from \eqref{e:Fun_Equ_LC_Fun_CInt}, while \eqref{e:F2p1_Od} derives from \eqref{e:Fun_Equ_LC_Fun_Sin} and \eqref{e:L_min_p}.
\end{proof}

\section{Series representation in the half-plane \texorpdfstring{$\sigma > 0$}{sigma > 0} for LC-functions with modulator vanishing at \texorpdfstring{0}{0}}
\label{s:Val_Ser_rep_LC_Fun_VM}

In this section we demonstrate that the series representation of an LC-function $L(s,f)$, as given by \eqref{e:LC_Fun_Ser}, is valid in the half-plane $\sigma>0$ provided that $p_f(0)=0$.

\begin{theorem}\label{t:Val_Ser_rep_LC_Fun_OM}
Let $L(s,f)$ be an LC-function whose modulator vanishes at $0$. The series defining $L(s,f)$ in \eqref{e:LC_Fun_Ser} converges absolutely for $\sigma>0$. The convergence is uniform in every compact subset of the half-plane $\sigma>0$. Moreover, since $n^{(-s,f)}$, for $n\geq n_f$, is a sequence of entire functions \cite[Theorem 4.2]{Lamgouni2024}, the series $\sum_{n=n_f}^{+\infty}n^{(-s, f)}$ represents an analytic function in the half-plane $\sigma>0$. Consequently, by the uniqueness of analytic continuation, \eqref{e:LC_Fun_Ser} holds for $\sigma>0$.
\end{theorem}

\begin{proof}
Given that $p_f(0)=0$, the first coefficient $P_{f,0}$ in the series \eqref{e:pf_Ser} is zero. Consequently, as stated in \eqref{e:Pow_szf2}, for all $n\geq n_f$, we have
\[
n^{(-s,f)}=\sum_{k=1}^{+\infty}\binom{-s}{k}P_{f,k}n^{-s-k}.
\]
According to \eqref{e:nf}, $r_f<n_f$. Let $r_0$ be an arbitrary real number such that $r_f<r_0<n_f$. Recall that $r_f=1/\rho_f$, where $\rho_f$ is the radius of convergence of the series $\sum_{k=0}^{+\infty}P_{f,k}z^k$. As $1/r_0<\rho_f$, the sequence $(P_{f,k}(1/r_0)^k)_{k \geq 1}$ is bounded by a constant $M>0$. From this, we deduce
\[
\abs{n^{(-s,f)}}\leq n^{-\sigma}\sum_{k=1}^{+\infty}\abs{\binom{-s}{k}}\abs{P_{f,k}\left(\frac{1}{r_0}\right)^{k}}\left(\frac{r_0}{n}\right)^{k}\leq Mn^{-\sigma}\sum_{k=1}^{+\infty}\abs{\binom{-s}{k}}\left(\frac{r_0}{n}\right)^{k}.
\]
Knowing that $\binom{-s}{k}=(-1)^k\prod_{l=0}^{k-1}(s+l)/k!$, we obtain
\[
\abs{\binom{-s}{k}}\leq\frac{\prod_{l=0}^{k-1}(\abs{s}+l)}{k!}=(-1)^k \binom{-\abs{s}}{k},
\]
which implies 
\[
\abs{n^{(-s,f)}}\leq Mn^{-\sigma}\sum_{k=1}^{+\infty}\binom{-\abs{s}}{k}\left(-\frac{r_0}{n}\right)^{k}.
\]
Since $0<r_0/n<1$, it follows that
\[
\abs{n^{(-s,f)}}\leq M n^{-\sigma}\left(\left(1-\frac{r_0}{n}\right)^{-\abs{s}}-1\right)=Mn^{-\sigma}\frac{n^{\abs{s}}-(n-r_0)^{\abs{s}}}{(n-r_0)^{\abs{s}}}.
\]
Applying the mean value theorem to the function $t\mapsto t^{\abs{s}}$ on $[n-r_0,n]$, we obtain
\[
n^{\abs{s}}-(n-r_0)^{\abs{s}}=\abs{s}r_0\lambda^{\abs{s}-1},
\]
where $n-r_0<\lambda<n$. Hence, \marginpar{\textcolor{red}{Corrected: replaced exponent $|s|$ with $|s|+1$ in inequality (5.1).}}
\begin{equation}\label{e:Ine_nsf}
\abs{n^{(-s,f)}}\leq\frac{Mr_0\abs{s}}{n^{\sigma+1}}\left(\frac{n}{n-r_0}\right)^{\abs{s}+1}\leq\frac{Mr_0\abs{s}}{n^{\sigma+1}}\left(\frac{n_f}{n_f-r_0}\right)^{\abs{s}+1}.
\end{equation}
This inequality ensures absolute convergence of the series $\sum_{n=n_f}^{+\infty} n^{(-s,f)}$ for $\sigma>0$.

Let $C$ be a compact subset of the half-plane $\sigma>0$. Then there exist constants $A>0$ and $\delta>0$ such that, for all $s\in C$, we have $\abs{s}<A$ and $\sigma\geq\delta$. From \eqref{e:Ine_nsf}, we obtain \marginpar{\textcolor{red}{Corrected: replaced exponent $A$ with $A+1$ in the adjacent inequality.}}
\[
\abs{n^{(-s,f)}}\leq\frac{AMr_0}{n^{\delta+1}}\left(\frac{n_f}{n_f-r_0}\right)^{A+1}.
\]
We conclude that the series $\sum_{n=n_f}^{+\infty}n^{(-s,f)}$ converges uniformly on $C$.
\end{proof}

\section{Applications}\label{s:App}

This section presents two examples of LC-functions with modulators of distinct parities, illustrating Theorem~\ref{t:Spe_Val_LC_fun_FC_fun}. We explicitly derive the values of these functions at even or odd positive integers, according to the parity of their modulators.

\subsection{Example of an LC-function with an even modulator}\label{ss:Exp1}

Consider the function $\varphi(t):=t\cos(wt)/(e^t-1)$, where $w$ is a complex parameter. First we keep $w$ real satisfying $0<\abs{w}<1$, and then extend the results to complex $w$ by analytic continuation. 

Let us first verify that the LC-function $L(s,\varphi)$ is  properly defined. By \eqref{e:pf_Ser}, the P-numbers associated to $\varphi$ are given by
\[
p_{\varphi}(t)=\cos(wt)=\sum_{n=0}^{+\infty}\frac{(-1)^n w^{2n}}{(2n)!} t^{2n}.
\]
That is, $P_{\varphi,2n+1}=0$ and $P_{\varphi,2n}=(iw)^{2n}$ for all integers $n\geq0$. Consequently, the series $\sum_{n=0}^{+\infty}P_{\varphi,n}z^{n}$ has a radius of convergence $\rho_{\varphi}=1/\abs{w}$, with $\rho_{\varphi}\to\infty$ as $w\to0$. Thus, $L(s,\varphi)$ is a well-defined LC-function with the even modulator $p_{\varphi}$.

\begin{figure}
\includegraphics{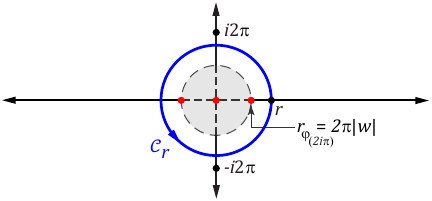}
\caption{The Hankel contour $\mu_{\varphi_{(2i\pi)}}$ reduces to a circle $\mathcal{C}_r$ with radius $r$ such that $2\pi\abs{w}<r<2\pi$.}
\label{f:Cer_FC_Fun_Phi}
\end{figure}

From \eqref{e:nf}, we have $n_{\varphi}=\lfloor\abs{w}\rfloor+1=1$ since by definition $r_{\varphi}=1/\rho_{\varphi}$. Now, as per \eqref{e:Pow_szf2}, for all $z\in\Omega_{\varphi}=\left\{\alpha\in\mathbb{C}\mid\alpha\notin\mathbb{R}_{\leq 0}\land\abs{\alpha}>\abs{w}\right\} $  and $ s\in\mathbb{C}$, the function $z^{(s,\varphi)}$ is written explicitly as
\[
z^{(s,\varphi)}=z^{s}\sum_{n=0}^{+\infty}\binom{s}{2n}\left(\frac{iw}{z}\right)^{2n}=z^{s}\frac{\left(1+\frac{iw}{z}\right)^{s}+\left(1-\frac{iw}{z}\right)^{s}}{2}.
\]
Hence, according to \eqref{e:LC_Fun_Ser}, the LC-function $ L(s,\varphi) $ is defined for $\sigma>1 $ by the series
\begin{equation}\label{e:LzPhi_Ser}
L(s,\varphi)=\frac{1}{2}\sum_{n=1}^{+\infty}\left(\frac{1}{(n+iw)^{s}}+\frac{1}{(n-iw)^{s}}\right).
\end{equation}
Alternatively, as stated in \eqref{e:LC_Fun_Int}, for $\sigma>1$, $L(s,\varphi)$ can be expressed by the integral
\begin{equation}\label{e:LzPhi_Int}
L(s,\varphi)=\frac{1}{\Gamma(s)}\int_0^{+\infty}\frac{t^{s-1}\cos(wt)}{e^t - 1}\,dt.
\end{equation}
On the other hand, referring to \eqref{e:L2p_Eve} in Theorem~\ref{t:Spe_Val_LC_fun_FC_fun}, we have for every non negative integer $p$
\[
L(2p,\varphi)=\frac{(-1)^{p}(2\pi)^{2p}}{4i\pi}\int_{\mu_{\varphi_{(2i\pi)}}}\frac{z^{(-2p,\varphi_{(2i\pi)})}}{e^{-z}-1}\,dz.
\]
According to \eqref{e:rfw_2Pirf}, the radius of the disk in Figure~\ref{f:Han_Con_FC_Fun_2iPi} is given by $r_{\varphi_{(2i\pi)}}=2\pi r_{\varphi}=2\pi\abs{w}$. The subsequent analysis aims to establish that the integrand of the integral above is a meromorphic function on $\mathbb{C}$ with all poles within the disk $\abs{z}\leq 2\pi\abs{w}$. This permits us to reduce the Hankel contour $\mu_{\varphi_{(2i\pi)}}$ to a circular path $\mathcal{C}_r$ centered at the origin, where the radius $r$ satisfies $2\pi\abs{w}<r<2\pi$ (see Figure~\ref{f:Cer_FC_Fun_Phi}). The final step involves applying Cauchy's residue theorem to derive a closed-form expression for $L(2p,\varphi)$.

To be more specific, From \eqref{e:Pow_szf2} and \eqref{e:Pfwn_wnPfn}, it follows that
\[
z^{(-2p,\varphi_{(2i\pi)})}=z^{-2p}\sum_{n=0}^{+\infty}\binom{-2p}{2n}(2i\pi)^{2n}P_{\varphi,2n}z^{-2n}=z^{-2p}\sum_{n=0}^{+\infty}\binom{-2p}{2n}\left(\frac{-2\pi w}{z}\right)^{2n}.
\]
On the Hankel path $\mu_{\varphi_{(2i\pi)}}$ (see Figure~\ref{f:Han_Con_FC_Fun_2iPi}), we have $\abs{z}>r_{\varphi_{(2i\pi)}}=2\pi\abs{w}$ as $z\in\Omega_{\varphi_{(2i\pi)}}$. Hence, $\abs{\frac{2\pi w}{z}}<1$. Therefore, we can write
\begin{align*}
z^{(-2p,\varphi_{(2i\pi)})}&=\frac{z^{-2p}}{2}\left[ \left(1-\frac{2\pi w}{z}\right)^{-2p}+\left(1+\frac{2\pi w}{z}\right)^{-2p}\right] \\
&=\frac{1}{2}\left[\frac{1}{(z-2\pi w)^{2p}}+\frac{1}{(z+2\pi w)^{2p}}\right].
\end{align*}
Then,
\[
L(2p,\varphi)=\frac{(-1)^{p}(2\pi)^{2p}}{8i\pi}\int_{\mathcal{C}_r}\frac{1}{e^{-z}-1}\left(\frac{1}{(z-2\pi w)^{2p}}+\frac{1}{(z+2\pi w)^{2p}}\right)\,dz.
\]

$-$ In the case where $p=0$, we immediately obtain
\[
L(0,\varphi)=\frac{1}{4i\pi}\int_{\mathcal{C}_r}\frac{dz}{e^{-z}-1}=\frac{1}{2}\underset{z=0}{\text{Res}}\frac{1}{e^{-z}-1}=-\frac{1}{2}.
\]

$-$ Assume now that $p\geq1$. In this case,
\begin{align*}
\begin{split}
L(2p,\varphi)=&\frac{(-1)^{p}(2\pi)^{2p}}{4}\left[\underset{z=0}{\text{Res}}\frac{1}{e^{-z}-1}\left(\frac{1}{(z-2\pi w)^{2p}}+\frac{1}{(z+2\pi w)^{2p}}\right)\right.\\
&\left.+\underset{z=2\pi w}{\text{Res}}\frac{1}{(e^{-z}-1)(z-2\pi w)^{2p}}+\underset{z=-2\pi w}{\text{Res}}\frac{1}{(e^{-z}-1)(z+2\pi w)^{2p}}\right].\\
=&\frac{(-1)^{p+1}}{2 w^{2p}}+\frac{(-1)^{p}(2\pi)^{2p}}{4(2p-1)!}\left(\lim_{2\pi w} h^{(2p-1)}(z)+\lim_{-2\pi w}h^{(2p-1)}(z)\right),
\end{split}
\end{align*}
where $ h $ is the function

\begin{equation}\label{e:h_Fun}
h(z):=\frac{1}{e^{-z}-1}.
\end{equation}

It is straightforward to verify by induction that the following formula for the $ n $-th derivative of the function $h$ holds for any integer $n\geq0$:
\begin{equation}\label{e:h_nth_Der}
h^{(n)}(z)=\frac{e^z A_n(e^z)}{(1-e^z)^{n+1}}.
\end{equation}
To establish this, it suffices to use the well-known recurrence relation satisfied by the Eulerian polynomials (see, e.g., Identity (50) in \cite[Page 214]{MR580155}, and Identity (7) in \cite[Page 34]{MR0272642})\footnote{The Eulerian polynomials $A_n(t)$ are represented by $A_{n,1}(t)$ in \cite[Equation (50)]{MR580155} and by ${}^1\!A_n(t)$ in \cite[Equation (7)]{MR0272642}.
}:
\[
(1+nt)A_n(t)+t(1-t)A'_n(t)=A_{n+1}(t).
\]
It follows then that
\[
L(2p,\varphi)=\frac{(-1)^{p+1}}{2 w^{2p}}+\frac{(-1)^{p}(2\pi)^{2p}}{4(2p-1)!}\left(\frac{e^{2\pi w}A_{2p-1}(e^{2\pi w})}{(1-e^{2\pi w})^{2p}}+ \frac{e^{-2\pi w}A_{2p-1}(e^{-2\pi w})}{(1-e^{-2\pi w})^{2p}}\right).
\]
The symmetry identity of the Eulerian polynomials (see, e.g., \cite[Equation (4), Page 30]{MR0272642}), 
\[
A_n(t)=t^{n-1}A_n(t^{-1}),
\] 
leads to
\begin{equation}\label{e:Eul_Pol_Sym}
\frac{e^{-z} A_n(e^{-z})}{(1-e^{-z})^{n+1}}=(-1)^{n+1}\frac{e^z A_n(e^z)}{(1-e^z)^{n+1}}.
\end{equation}
As a result, we obtain
\[
L(2p,\varphi)=\frac{(2i\pi)^{2p}}{2(2p-1)!}\frac{e^{2\pi w}A_{2p-1}(e^{2\pi w})}{(1-e^{2\pi w})^{2p}}-\frac{(-1)^p}{2w^{2p}}.
\]
By introducing the function
\begin{equation}\label{e:K_Fun}
K_m(z):=\frac{(2i)^{m}e^{2z}A_{m-1}(e^{2z})}{(1-e^{2z})^{m}},
\end{equation}
where $m$ is a positive integer, the last identity can be expressed as \marginpar{\textcolor{red}{Corrected: replaced $\pi^{2^p}$ with $\pi^{2p}$ in the adjacent equality.}}
\begin{equation}\label{e:L2pPhi_Exp}
L(2p,\varphi)=\frac{\pi^{2p}}{2(2p-1)!}K_{2p}(\pi w)-\frac{(-1)^{p}}{2w^{2p}}.
\end{equation}

Given the two representations of $L(2p,\varphi)$, one as a series in \eqref{e:LzPhi_Ser} and the other as an integral in \eqref{e:LzPhi_Int}, we derive the following identities: \marginpar{\textcolor{red}{Corrected: replaced $\pi^{2^p}$ with $\pi^{2p}$ in the adjacent equality.}}
\begin{align}\label{e:K2p_Piw_Ser}
\frac{\pi^{2p}}{(2p-1)!} K_{2p}(\pi w)&=\frac{1}{(iw)^{2p}}+\sum_{n=1}^{+\infty}\left(\frac{1}{(n + iw)^{2p}}+\frac{1}{(n-iw)^{2p}}\right)\\
&=\frac{1}{(iw)^{2p}}+\zeta(2p,1+iw)+\zeta(2p,1-iw),\nonumber
\end{align}
and
\begin{equation}\label{e:Int_Cos_Exp}
\int_0^{+\infty}\frac{t^{2p-1}\cos(wt)}{e^t-1}\,dt=\frac{(-1)^p}{2}\left(\frac{(2\pi)^{2p}e^{2\pi w}A_{2p-1}(e^{2\pi w})}{(1-e^{2\pi w})^{2p}}- \frac{(2p-1)!}{w^{2p}}\right).
\end{equation}
Both members in \eqref{e:K2p_Piw_Ser} agree on $(-1,1)\setminus\{0\}$, implying by the identity theorem for holomorphic functions (see, e.g., \cite[Chapter 8, Page 228]{MR1084167} and \cite[Page 241]{MR0716497}) that they agree for all complex $w\in\mathbb{C}\setminus\left\{ik\mid k\in\mathbb{Z}\right\}$. By making the substitution $iw$ in place of $w$, we reformulate \eqref{e:K2p_Piw_Ser} for $w\in\mathbb{C}\setminus\mathbb{Z}$ as follows \marginpar{\textcolor{red}{Corrected: replaced $\pi^{2^p}$ with $\pi^{2p}$ in the adjacent equality.}}
\begin{align}\label{e:K2p_iPiw_Ser}
\frac{\pi^{2p}}{(2p-1)!} K_{2p}(i\pi w)&=\frac{1}{w^{2p}}+\sum_{n=1}^{+\infty}\left(\frac{1}{(w-n)^{2p}}+\frac{1}{(w+n)^{2p}}\right)\\
&=\frac{1}{w^{2p}}+\zeta(2p,1-w)+\zeta(2p,1+w). \nonumber
\end{align}
Once again, it follows by the analytic continuation that \eqref{e:Int_Cos_Exp} holds for all complex $w$ satisfying $\abs{\Im(w)}<1$. The singularity at $0$ of the right member is removable, and for $p\geq1$ we have the limit
\[
\lim_{w\to0}\frac{(-1)^p}{2}\left(\frac{(2\pi)^{2p}e^{2\pi w}A_{2p-1}(e^{2\pi w})}{(2p-1)!(1-e^{2\pi w})^{2p}}-\frac{1}{w^{2p}}\right)=\zeta(2p).
\]

 \subsection{Example of an LC-function with an odd modulator}\label{ss:Exp2}
The analysis of this example follows the same procedure as the first example, and is presented here in outline form only. Consider the LC-function $L(s,\phi)$ associated to the function $\phi(t):=t\sinh(wt)/(e^t-1)$, where $w$ is a complex parameter initially assumed to be real such that $0<\abs{w}<1$. The modulator of $L(s,\phi)$ and the P-numbers associated to $\phi$ are given by
\[
p_{\phi}(t)=\sinh(wt)=\sum_{n=0}^{+\infty}\frac{w^{2n+1}}{(2n+1)!}t^{2n+1}.
\]
From \eqref{e:nf}, we deduce that $n_{\phi}=1$, as $r_{\phi}=\abs{w}<1$. Furthermore, according to \eqref{e:Pow_szf2}, for all $z\in\Omega_{\phi}$ and $s\in\mathbb{C}$, since $\abs{z}>\abs{w}$, we have 
\[
z^{(s,\phi)}=z^{s}\sum_{n=0}^{+\infty}\binom{s}{2n+1}\left(\frac{w}{z}\right)^{2n+1}=z^{s}\frac{\left(1+\frac{w}{z}\right)^{s}-\left(1-\frac{w}{z}\right)^{s}}{2}.
\]
Therefore, by \eqref{e:LC_Fun_Ser} and \eqref{e:LC_Fun_Int}, $L(s,\phi)$ is given for $\sigma>1$ by the two equivalent expressions
\begin{equation}\label{e:LzPsi_Ser_Int}
L(s,\phi)=\frac{1}{2}\sum_{n=1}^{+\infty}\left(\frac{1}{(n+w)^{s}}-\frac{1}{(n-w)^{s}}\right)=\frac{-1}{\Gamma(s)}\int_0^{+\infty}\frac{t^{s-1}\sinh(wt)}{e^t-1}\,dt.
\end{equation}
Due to Remark~\ref{r:Ana_Con_LC_Fun}.(ii), the LC-function $L(s,\phi)$ is entire because $p_{\phi}$ is odd. Furthermore, the series representation above holds for $\sigma>0$, as demonstrated in Theorem~\ref{t:Val_Ser_rep_LC_Fun_OM}.

Now, from \eqref{e:L2p1_Od} in Theorem~\ref{t:Spe_Val_LC_fun_FC_fun}, we obtain for every non negative integer $p$
\begin{align*}
L(2p+1,\phi)&=\frac{(-1)^p(2\pi)^{2p+1}}{4\pi}\int_{\mu_{\phi_{(2i\pi)}}}\frac{z^{(-2p-1,\phi_{(2i\pi)})}}{e^{-z}-1}\,dz \\
&=\frac{1}{2i\pi}\int_{\mathcal{C}_r}\frac{(2i\pi)^{2p+1}}{4(e^{-z} - 1)}\left( \frac{1}{(z + 2i\pi w)^{2p+1}}-\frac{1}{(z-2i\pi w)^{2p+1}}\right) \,dz\\
&=-\frac{1}{2 w^{2p+1}}+\frac{(2i\pi)^{2p+1}}{4(2p)!}\left(\lim_{-2i\pi w}h^{(2p)}(z)-\lim_{2i\pi w}h^{(2p)}(z) \right).
\end{align*} 
Here, $h$ is the function considered in \eqref{e:h_Fun}.

$-$ For $p=0$, we obtain
\begin{equation}\label{e:L1Psi}
L(1,\phi)=-\frac{1}{2w}+\frac{\pi}{2}\cot(\pi w).
\end{equation}
By equating the series representation of $ L(1,\phi) $ given in \eqref{e:LzPsi_Ser_Int}, with its explicit value provided in \eqref{e:L1Psi}, we find Euler's formula for the partial fraction expansion of $\cot(\pi w)$, as stated in \eqref{e:ParFraExpCot}.

$-$ For $p\geq1$ we have, from \eqref{e:h_nth_Der} and \eqref{e:Eul_Pol_Sym},
\marginpar{\textcolor{red}{Shortened the equation below from three lines to two.}}
\begin{align*}
L(2p+1,\phi) &= -\frac{1}{2 w^{2p+1}}+\frac{(2i\pi)^{2p+1}}{4(2p)!}\left(\frac{e^{-2i\pi w}A_{2p}(e^{-2i\pi w})}{(1-e^{-2i\pi w})^{2p+1}}-\frac{e^{2i\pi w}A_{2p}(e^{2i\pi w})}{(1-e^{2i\pi w})^{2p+1}}\right)\\
&=-\frac{(2i\pi)^{2p+1}}{2(2p)!}\frac{e^{2i\pi w} A_{2p}(e^{2i\pi w})}{(1-e^{2i\pi w})^{2p+1}}-\frac{1}{2 w^{2p+1}}.
\end{align*}
According to the definition of the function $K_m(z)$ in \eqref{e:K_Fun}, we write
\begin{equation}\label{e:L2p1Psi_Exp}
L(2p+1,\phi)=-\frac{\pi^{2p+1}}{2(2p)!}K_{2p+1}(i\pi w)-\frac{1}{2w^{2p+1}}.
\end{equation}
Given the two representations of $L(2p+1,\phi)$ in \eqref{e:LzPsi_Ser_Int}, we derive the following results:
\begin{align}\label{e:Int_Sinh_Exp}
\int_0^{+\infty}\frac{t^{2p}\sinh(wt)}{e^t-1}\,dt&=\frac{\pi^{2p+1}}{2}K_{2p+1}(i\pi w)+\frac{(2p)!}{2w^{2p+1}}\\
&=\frac{1}{2}\left(\frac{(2i\pi)^{2p+1}e^{2i\pi w} A_{2p}(e^{2i\pi w})}{(1-e^{2i\pi w})^{2p+1}}+\frac{(2p)!}{w^{2p+1}}\right)\nonumber,
\end{align}
and
\begin{equation}\label{e:K2p1_iPiw_Ser}
\frac{\pi^{2p+1}}{(2p)!} K_{2p+1}(i\pi w)=-\frac{1}{w^{2p+1}}-\sum_{n=1}^{+\infty}\left(\frac{1}{(w-n)^{2p+1}}+\frac{1}{(w+n)^{2p+1}}\right).
\end{equation}
These identities can be extended to complex $w$ by analytic continuation; \eqref{e:Int_Sinh_Exp} holds for all complex $w$ satisfying $\abs{\Re(w)}<1$, with 
\[
\lim_{w\to0}\frac{(2i\pi)^{2p+1}e^{2i\pi w}A_{2p}(e^{2i\pi w})}{(1-e^{2i\pi w})^{2p+1}}+\frac{(2p)!}{w^{2p+1}}=0,
\]
and \eqref{e:K2p1_iPiw_Ser} holds for all $w\in\mathbb{C}\setminus\mathbb{Z}$. 

The identities \eqref{e:K2p1_iPiw_Ser} and \eqref{e:K2p_iPiw_Ser} can be unified into a single expression, valid for all $w\in\mathbb{C}\setminus\mathbb{Z}$ and for all integer $p\geq2$:
\begin{equation}\label{e:_Kp_iPiw_Ser}
\frac{(-\pi)^{p}}{(p-1)!} K_p(i\pi w)=\frac{1}{w^{p}}+\sum_{n=1}^{+\infty}\left(\frac{1}{(w-n)^{p}}+\frac{1}{(w+n)^{p}}\right).
\end{equation}

Extending our analysis, we express the function $K_p (i\pi w)$ in terms of $\cos(\pi w)$, $\sin(\pi w)$, and Eulerian numbers. Subsequently, we establish the identity \[\cot^{(p)}(\pi w)=-K_{p+1}(i\pi w),\] thus concluding the example.

According to equation \eqref{e:K_Fun}, we have
\[
K_p(i\pi w)=\frac{(2i)^{p}e^{i(2-p)\pi w}A_{p-1}(e^{2i\pi w})}{(e^{-i\pi w}-e^{i\pi w})^p}=\frac{(-1)^{p}e^{i(2-p)\pi w}}{\sin^p(\pi w)}\sum_{k=0}^{p-2}\eulerianNum{p-1}{k}e^{2ik\pi w},
\]
where $\eulerianNum{\cdot}{\cdot}$ denotes the Eulerian numbers, generated by the Eulerian polynomials (see, e.g., \cite[Equation (1.8)]{MR3408615} and \cite[Page 264]{MR2744266}). We distinguish two cases:

\textsc{First case}. $p$ is even; $p=2N$ with $N\geq1$. We formulate 
\begin{align*}
K_{2N}(i\pi w)&=\frac{1}{\sin^{2N}(\pi w)}\left(\sum_{k=0}^{N-2}\eulerianNum{2N-1}{k}e^{-2i(N-k-1)\pi w}\right. \\
&\quad\left.+\eulerianNum{2N-1}{N-1}+\sum_{k=N}^{2N-2}\eulerianNum{2N-1}{k}e^{2i(k-N+1)\pi w}\right) \\
&=\frac{1}{\sin^{2N}(\pi w)}\left(\sum_{k=1}^{N-1}\eulerianNum{2N-1}{N-k-1}e^{-2ik\pi w}\right. \\
&\quad\left.+\eulerianNum{2N-1}{N-1}+\sum_{k=1}^{N-1}\eulerianNum{2N-1}{N+k-1}e^{2ik\pi w}\right).
\end{align*}
Therefore using the identity $\eulerianNum{n}{k}=\eulerianNum{n}{n-k-1}$ we conclude the final form
\[
K_{2N}(i\pi w)=\frac{\eulerianNum{2N-1}{N-1}+2\sum_{k=1}^{N-1}\eulerianNum{2N-1}{N-k-1}\cos(2k\pi w)}{\sin^{2N}(\pi w)}.
\]

\textsc{Second case}. $p$ is odd; $p=2N+1$ with $N\geq1$. By the same calculation as in the first case we obtain
\[
K_{2N+1}(i\pi w)=-\frac{2\sum_{k=1}^{N}\eulerianNum{2N}{N-k}\cos((2k-1)\pi w)}{\sin^{2N+1}(\pi w)}.
\]
From \eqref{e:_Kp_iPiw_Ser}, for $N\in\mathbb{Z}_{\geq 1}$ and $w\in\mathbb{C}\setminus\mathbb{Z}$ we obtain the explicit formulas:
\begin{multline}
\frac{\eulerianNum{2N-1}{N-1}+2\sum_{k=1}^{N-1}\eulerianNum{2N-1}{N-k-1}\cos(2k\pi w)}{(2N-1)!\sin^{2N}(\pi w)}\pi^{2N}\\
=\frac{1}{w^{2N}}+\sum_{n=1}^{+\infty}\left(\frac{1}{(w-n)^{2N}}+\frac{1}{(w+n)^{2N}}\right),
\end{multline}
\begin{multline}
\frac{2\sum_{k=1}^{N}\eulerianNum{2N}{N-k}\cos((2k-1)\pi w)}{(2N)!\sin^{2N+1}(\pi w)}\pi^{2N+1}\\
=\frac{1}{w^{2N+1}}+\sum_{n=1}^{+\infty}\left(\frac{1}{(w-n)^{2N+1}}+\frac{1}{(w+n)^{2N+1}}\right).
\end{multline}

Now, by differentiating Euler's formula \eqref{e:ParFraExpCot} to the order $p\geq1$, then matching the result with \eqref{e:_Kp_iPiw_Ser}, we deduce that 
\begin{align}\label{e:DerCotOrd_p}
\cot^{(p)}(\pi w)&=\frac{(-1)^{p}p!}{\pi^{p+1}}\left(\frac{1}{w^{p+1}}+\sum_{n=1}^{+\infty}\left(\frac{1}{(w-n)^{p+1}}+ \frac{1}{(w+n)^{p+1}}\right)\right)\\
&=-K_{p+1}(i\pi w)=-\frac{(2i)^{p+1}e^{2i\pi w}A_p(e^{2i\pi w})}{(1-e^{2i\pi w})^{p+1}}.\nonumber
\end{align}
This yields expressions for the higher derivatives of $\cot(w)$ in terms of $\cos$, $\sin$, and Eulerian numbers, valid for $w\in\mathbb{C}\setminus\pi\mathbb{Z}$:
\begin{align}
\cot^{(2N)}(w)&=\frac{2\sum_{k=1}^{N}\eulerianNum{2N}{N-k}\cos((2k-1)w)}{\sin^{2N+1}(w)},\\
\cot^{(2N-1)}(w)&=-\frac{\eulerianNum{2N-1}{N-1}+2\sum_{k=1}^{N-1}\eulerianNum{2N-1}{N-k-1}\cos(2kw)}{\sin^{2N}(w)}.
\end{align}

\bibliographystyle{aomplain}
\bibliography{AMS-bibliography}

@Article{Lamgouni2024,
author={Lamgouni, Lahcen},
title={C-polynomials and {LC}-functions: towards a generalization of the {H}urwitz zeta function},
journal={The Ramanujan Journal},
year={2024},
month={Aug},
day={16},
issn={1572-9303},
doi={10.1007/s11139-024-00919-1},
url={https://doi.org/10.1007/s11139-024-00919-1}
}

@incollection {MR1811017,
     AUTHOR = {Garunk{\v{s}}tis, R. and Laurin{\v{c}}ikas, A.},
     TITLE = {The {L}erch zeta-function},
      NOTE = {Analytical methods of analysis and differential equations
              (Minsk, 1999)},
   JOURNAL = {Integral Transform. Spec. Funct.},
  FJOURNAL = {Integral Transforms and Special Functions. An International
              Journal},
    VOLUME = {10},
      YEAR = {2000},
    NUMBER = {3-4},
     PAGES = {211--226},
      ISSN = {1065-2469,1476-8291},
   MRCLASS = {11M35},
  MRNUMBER = {1811017},
       DOI = {10.1080/10652460008819287},
       URL = {https://doi.org/10.1080/10652460008819287},
}

@article {MR43843,
    AUTHOR = {Apostol, T. M.},
     TITLE = {On the {L}erch zeta function},
   JOURNAL = {Pacific J. Math.},
  FJOURNAL = {Pacific Journal of Mathematics},
    VOLUME = {1},
      YEAR = {1951},
     PAGES = {161--167},
      ISSN = {0030-8730,1945-5844},
   MRCLASS = {10.0X},
  MRNUMBER = {43843},
MRREVIEWER = {N.\ G.\ de Bruijn},
       URL = {http://projecteuclid.org/euclid.pjm/1103052188},
}

@book {MR1697859,
    AUTHOR = {Neukirch, J\"urgen},
     TITLE = {Algebraic number theory},
    SERIES = {Grundlehren der mathematischen Wissenschaften [Fundamental
              Principles of Mathematical Sciences]},
    VOLUME = {322},
      NOTE = {Translated from the 1992 German original and with a note by
              Norbert Schappacher,
              With a foreword by G. Harder},
 PUBLISHER = {Springer-Verlag, Berlin},
      YEAR = {1999},
     PAGES = {xviii+571},
      ISBN = {3-540-65399-6},
   MRCLASS = {11Rxx (11-02 11S15 11S31 14C40)},
  MRNUMBER = {1697859},
MRREVIEWER = {Cornelius\ Greither},
       DOI = {10.1007/978-3-662-03983-0},
       URL = {https://doi.org/10.1007/978-3-662-03983-0},
}

@book {MR0606931,
    AUTHOR = {Davenport, Harold},
     TITLE = {Multiplicative number theory},
    SERIES = {Graduate Texts in Mathematics},
    VOLUME = {74},
   EDITION = {Second},
      NOTE = {Revised by Hugh L. Montgomery},
 PUBLISHER = {Springer-Verlag, New York-Berlin},
      YEAR = {1980},
     PAGES = {xiii+177},
      ISBN = {0-387-90533-2},
   MRCLASS = {10-01 (10-02 10Hxx)},
  MRNUMBER = {606931},
MRREVIEWER = {H.-E.\ Richert},
}

@article {MR2860694,
    AUTHOR = {Alkan, Emre},
     TITLE = {Values of {D}irichlet {$L$}-functions, {G}auss sums and
              trigonometric sums},
   JOURNAL = {Ramanujan J.},
  FJOURNAL = {Ramanujan Journal. An International Journal Devoted to the
              Areas of Mathematics Influenced by Ramanujan},
    VOLUME = {26},
      YEAR = {2011},
    NUMBER = {3},
     PAGES = {375--398},
      ISSN = {1382-4090,1572-9303},
   MRCLASS = {11M06 (11L03 11L05)},
  MRNUMBER = {2860694},
MRREVIEWER = {Moubariz\ Z.\ Garaev},
       DOI = {10.1007/s11139-010-9292-8},
       URL = {https://doi.org/10.1007/s11139-010-9292-8},
}

@article {MR1874365,
    AUTHOR = {Louboutin, St\'ephane},
     TITLE = {The mean value of {$|L(k,\chi)|^2$} at positive rational
              integers {$k\geq1$}},
   JOURNAL = {Colloq. Math.},
  FJOURNAL = {Colloquium Mathematicum},
    VOLUME = {90},
      YEAR = {2001},
    NUMBER = {1},
     PAGES = {69--76},
      ISSN = {0010-1354,1730-6302},
   MRCLASS = {11M06 (11M20)},
  MRNUMBER = {1874365},
MRREVIEWER = {J\"orn\ Steuding},
       DOI = {10.4064/cm90-1-6},
       URL = {https://doi.org/10.4064/cm90-1-6},
}

@article {MR1222534,
    AUTHOR = {Louboutin, St\'ephane},
     TITLE = {Quelques formules exactes pour des moyennes de fonctions {$L$} de {D}irichlet},
   JOURNAL = {Canad. Math. Bull.},
  FJOURNAL = {Canadian Mathematical Bulletin. Bulletin Canadien de
              Math\'ematiques},
    VOLUME = {36},
      YEAR = {1993},
    NUMBER = {2},
     PAGES = {190--196},
      ISSN = {0008-4395,1496-4287},
   MRCLASS = {11M20 (11R18)},
  MRNUMBER = {1222534},
MRREVIEWER = {D.\ R.\ Heath-Brown},
       DOI = {10.4153/CMB-1993-028-8},
       URL = {https://doi.org/10.4153/CMB-1993-028-8},
}

@book {MR2376618,
    AUTHOR = {Murty, M. Ram},
     TITLE = {Problems in analytic number theory},
    SERIES = {Graduate Texts in Mathematics},
    VOLUME = {206},
   EDITION = {Second},
      NOTE = {Readings in Mathematics},
 PUBLISHER = {Springer, New York},
      YEAR = {2008},
     PAGES = {xxii+502},
      ISBN = {978-0-387-72349-5},
   MRCLASS = {11-01 (11M06 11N05)},
  MRNUMBER = {2376618},
}

@article {MR0285499,
    AUTHOR = {Apostol, Tom M.},
     TITLE = {Dirichlet {$L$}-functions and primitive characters},
   JOURNAL = {Proc. Amer. Math. Soc.},
  FJOURNAL = {Proceedings of the American Mathematical Society},
    VOLUME = {31},
      YEAR = {1972},
     PAGES = {384--386},
      ISSN = {0002-9939,1088-6826},
   MRCLASS = {10.41},
  MRNUMBER = {285499},
MRREVIEWER = {K.\ Thanigasalam},
       DOI = {10.2307/2037537},
       URL = {https://doi.org/10.2307/2037537},
}

@book {MR1084167,
    AUTHOR = {Remmert, Reinhold},
     TITLE = {Theory of complex functions},
    SERIES = {Graduate Texts in Mathematics},
    VOLUME = {122},
   EDITION = {German},
      NOTE = {Readings in Mathematics},
 PUBLISHER = {Springer-Verlag, New York},
      YEAR = {1991},
     PAGES = {xx+453},
      ISBN = {0-387-97195-5},
   MRCLASS = {30-01},
  MRNUMBER = {1084167},
MRREVIEWER = {M.\ H.\ Heins},
       DOI = {10.1007/978-1-4612-0939-3},
       URL = {https://doi.org/10.1007/978-1-4612-0939-3},
}

@book {MR0434929,
    AUTHOR = {Apostol, Tom M.},
     TITLE = {Introduction to analytic number theory},
    SERIES = {Undergraduate Texts in Mathematics},
 PUBLISHER = {Springer-Verlag, New York-Heidelberg},
      YEAR = {1976},
     PAGES = {xii+338},
   MRCLASS = {10-01 (10AXX 10HXX)},
  MRNUMBER = {434929},
MRREVIEWER = {E.\ Grosswald},
}

@book {MR3408615,
    AUTHOR = {Petersen, T. Kyle},
     TITLE = {Eulerian numbers},
    SERIES = {Birkh\"auser Advanced Texts: Basler Lehrb\"ucher.
              [Birkh\"auser Advanced Texts: Basel Textbooks]},
      NOTE = {With a foreword by Richard Stanley},
 PUBLISHER = {Birkh\"auser/Springer, New York},
      YEAR = {2015},
     PAGES = {xviii+456},
      ISBN = {978-1-4939-3090-6; 978-1-4939-3091-3},
   MRCLASS = {05-02 (05A15 05Exx 06A07 11B65 11B75 20F55)},
  MRNUMBER = {3408615},
MRREVIEWER = {Damir\ Yeliussizov},
       DOI = {10.1007/978-1-4939-3091-3},
       URL = {https://doi.org/10.1007/978-1-4939-3091-3},
}

@book {MR1216135,
    AUTHOR = {Hida, Haruzo},
     TITLE = {Elementary theory of {$L$}-functions and {E}isenstein series},
    SERIES = {London Mathematical Society Student Texts},
    VOLUME = {26},
 PUBLISHER = {Cambridge University Press, Cambridge},
      YEAR = {1993},
     PAGES = {xii+386},
      ISBN = {0-521-43411-4; 0-521-43569-2},
   MRCLASS = {11Fxx (11Mxx 11R42)},
  MRNUMBER = {1216135},
MRREVIEWER = {Jacques\ Tilouine},
       DOI = {10.1017/CBO9780511623691},
       URL = {https://doi.org/10.1017/CBO9780511623691},
}

@book {MR2378655,
    AUTHOR = {Montgomery, Hugh L. and Vaughan, Robert C.},
     TITLE = {Multiplicative number theory. {I}. {C}lassical theory},
    SERIES = {Cambridge Studies in Advanced Mathematics},
    VOLUME = {97},
 PUBLISHER = {Cambridge University Press, Cambridge},
      YEAR = {2007},
     PAGES = {xviii+552},
      ISBN = {978-0-521-84903-6; 0-521-84903-9},
   MRCLASS = {11-02 (11-01 11M06 11M26 11M45 11N05 11N37)},
  MRNUMBER = {2378655},
MRREVIEWER = {Wolfgang\ Schwarz},
}

@incollection {MR2744266,
    AUTHOR = {Foata, Dominique},
     TITLE = {Eulerian polynomials: from {E}uler's time to the present},
 BOOKTITLE = {The legacy of {A}lladi {R}amakrishnan in the mathematical
              sciences},
     PAGES = {253--273},
 PUBLISHER = {Springer, New York},
      YEAR = {2010},
      ISBN = {978-1-4419-6262-1},
   MRCLASS = {01A50 (05A15 05A30 33B10)},
  MRNUMBER = {2744266},
MRREVIEWER = {E.\ Keith\ Lloyd},
       DOI = {10.1007/978-1-4419-6263-8\_15},
       URL = {https://doi.org/10.1007/978-1-4419-6263-8_15},
}

@book {MR2569612,
    AUTHOR = {Aigner, Martin and Ziegler, G\"unter M.},
     TITLE = {Proofs from {T}he {B}ook},
   EDITION = {Fourth},
 PUBLISHER = {Springer-Verlag, Berlin},
      YEAR = {2010},
     PAGES = {viii+274},
      ISBN = {978-3-642-00855-9},
   MRCLASS = {00A05},
  MRNUMBER = {2569612},
       DOI = {10.1007/978-3-642-00856-6},
       URL = {https://doi.org/10.1007/978-3-642-00856-6},
}

@book {MR0716497,
    AUTHOR = {Kaup, Ludger and Kaup, Burchard},
     TITLE = {Holomorphic functions of several variables},
    SERIES = {De Gruyter Studies in Mathematics},
    VOLUME = {3},
      NOTE = {An introduction to the fundamental theory,
              With the assistance of Gottfried Barthel,
              Translated from the German by Michael Bridgland},
 PUBLISHER = {Walter de Gruyter \& Co., Berlin},
      YEAR = {1983},
     PAGES = {xv+349},
      ISBN = {3-11-004150-2},
   MRCLASS = {32-01},
  MRNUMBER = {716497},
MRREVIEWER = {R.\ Michael\ Range},
       DOI = {10.1515/9783110838350},
       URL = {https://doi.org/10.1515/9783110838350},
}

@book {MR580155,
    AUTHOR = {Riordan, John},
     TITLE = {An introduction to combinatorial analysis},
      NOTE = {Reprint of the 1958 edition},
 PUBLISHER = {Princeton University Press, Princeton, NJ},
      YEAR = {1980},
     PAGES = {xii+244},
      ISBN = {0-691-02365-4},
   MRCLASS = {05-01},
  MRNUMBER = {580155},
}

@book {MR0272642,
    AUTHOR = {Foata, Dominique and Sch\"utzenberger, Marcel-P.},
     TITLE = {Th\'eorie g\'eom\'etrique des polyn\^omes eul\'eriens},
    SERIES = {Lecture Notes in Mathematics},
    VOLUME = {Vol. 138},
 PUBLISHER = {Springer-Verlag, Berlin-New York},
      YEAR = {1970},
     PAGES = {v+94},
   MRCLASS = {05.10},
  MRNUMBER = {272642},
MRREVIEWER = {D.\ P.\ Roselle},
}

@book {MR1841793,
    AUTHOR = {Euler, Leonhard},
     TITLE = {Introductio in analysin infinitorum. {T}omus primus},
      NOTE = {Reprint of the 1748 original},
 PUBLISHER = {Sociedad Andaluza de Educaci\'on Matem\'atica ``Thales'',
              Seville; Real Sociedad Matem\'atica Espa\~nola, Madrid},
      YEAR = {2000},
     PAGES = {xvi+320},
      ISBN = {84-923760-2-3},
   MRCLASS = {01A50 (01A75)},
  MRNUMBER = {1841793},
}

@manual{orcidlink,
  title = {The orcidlink package},
  subtitle = {Insert hyperlinked ORCiD logo},
  author = {Stein, Leo},
  url = {https://ctan.org/pkg/orcidlink},
  urldate = {2024-06-08}, 
  yar = 2020,
  version = {1.0.5},
}

\end{document}